\documentclass[draft]{publmathdeb}
\usepackage{amsmath,amsfonts,amssymb,hyperref, vhistory}
\author{Akhilesh Yadav}
\author{Kiran Meena}
\address{Department of Mathematics\\ Institute of Science\\
	Banaras Hindu University\\
	Varanasi-221005, India}
\email{akhilesha68@gmail.com (A. Yadav)}
\email{kirankapishmeena@gmail.com (K. Meena)}
\title[Riemannian maps whose base manifolds admit a Ricci soliton]
{Riemannian maps whose base manifolds admit a Ricci soliton}
\thanks{The second author gratefully acknowledges the financial support provided by the Council of	Scientific and Industrial Research – Human Resource Development Group (CSIR-HRDG), New Delhi, India (File No.: 09/013(0887)/2019-EMR-I)}
\begin{document}
	\begin{abstract}
		In this paper, we study Riemannian maps whose base manifolds admit a Ricci soliton and give a non-trivial example of such a Riemannian map. First, we find Riemannian curvature tensor for the base manifolds of Riemannian map $F$. Further, we obtain the Ricci tensor and calculate the scalar curvature of the base manifold. Moreover, we obtain necessary conditions for the leaves of $rangeF_\ast$ to be Ricci soliton, almost Ricci soliton, and Einstein. We also obtain necessary conditions for the leaves of $(rangeF_\ast)^\bot$ to be Ricci soliton and Einstein. Also, we calculate the scalar curvatures of $rangeF_\ast$ and $(rangeF_\ast)^\bot$ by using Ricci soliton. Finally, we study the harmonicity and biharmonicity of such a  Riemannian map. We obtain a necessary and sufficient condition for such a Riemannian map between Riemannian manifolds to be harmonic. We also obtain necessary and sufficient conditions for a Riemannian map from a Riemannian manifold to a space form that admits Ricci soliton to be harmonic and biharmonic.
	\end{abstract}
	
	\maketitle
	\theoremstyle{plain}
	\theoremstyle{definition}
	\theoremstyle{remark}
	\theoremstyle{plain}
	\newtheorem{theorem}{Theorem}[section]
	\newtheorem{proposition}[theorem]{Proposition}
	\newtheorem{lemma}[theorem]{Lemma}
	\newtheorem{corollary}[theorem]{Corollary}
	\newtheorem{state}{State}[section]
	\theoremstyle{definition}
	\newtheorem{definition}{Definition}[section]
	\theoremstyle{remark}
	\newtheorem{remark}{Remark}[section]
	\newtheorem{example}{Example}
	\theoremstyle{plain}
	\newtheorem*{hth}{Henk's Theorem}
	\section{Introduction}
	In 1992, Fischer introduced the notion of Riemannian map between Riemannian manifolds in \cite{Fischer} as a generalization of the notion of an isometric immersion and Riemannian submersion. The geometry of Riemannian submersions has been discussed in \cite{Falcitelli}. We note that a remarkable property of Riemannian maps is that a Riemannian map satisfies the generalized eikonal equation $\| F_\ast \|^2 = rank F$, which is a bridge between geometric optics and physical optics \cite{Fischer}. The eikonal equation of geometrical optics was solved by using Cauchy's method of characteristics. In \cite{Fischer} Fischer also proposed an approach to building a quantum model and he pointed out the success of such a program of building a quantum model of nature using Riemannian maps. It provide an interesting relationship between Riemannian maps, harmonic maps, and Lagrangian field theory on the mathematical side, and Maxwell's equation, Shr\"{o}dinger's equation and their proposed generalization on the physical side.
	
	In \cite{Sahin2015}, B. \c{S}ahin developed certain geometric structures along a Riemannian map to investigate the geometry of such a map. He constructed Gauss-Weingarten formulas and obtained Gauss, Codazzi and Ricci equations for Riemannian map by using the second fundamental form and suitable linear connections.
	
	On the other hand, in 1988, the notion of Ricci soliton was introduced by Hamilton in \cite{Hamilton}. A Ricci soliton is a natural generalization of an Einstein metric. A Riemannian manifold $(N,g_N)$ is called a Ricci soliton if there exists a smooth vector field $\xi$ (called potential vector field) on $N$ such that 
	\begin{equation}\label{eqn1.1} 
		\frac{1}{2} (L_{\xi} g_N)(X_1,Y_1) + Ric (X_1,Y_1) + \lambda g_N(X_1,Y_1) = 0,
	\end{equation}
	where $L_{\xi}$ is the Lie derivative of the metric tensor of $g_N$ with respect to $\xi$, $Ric$ is the Ricci tensor of $(N,g_N)$, $\lambda$ is a constant and $X_1$, $Y_1$ are arbitrary vector fields on $N$. We shall denote a Ricci soliton by $(N,g_{N},\xi,\lambda)$. The Ricci soliton $(N,g_N,\xi,\lambda)$ is said to be shrinking, steady or expanding accordingly as $\lambda<0$, $\lambda=0$ or $\lambda>0$, respectively. It is obvious that a trivial Ricci soliton is an Einstein manifold with $\xi$ zero or killing, that is, Lie derivative of metric tensor $g_N$ with respect to $\xi$ vanishes. Hamilton showed that the self-similar solutions of Ricci flow are Ricci solitons. The Ricci soliton $(N,g_N,\xi,\lambda)$ is said to be a gradient Ricci soliton if the potential vector field $\xi$ is the gradient of some smooth function $f$ on $N$, which is denoted by $(N,g_N,f,\lambda)$. A non-killing tangent vector field $\xi$ on a Riemannian manifold $(N,g_N)$ is called conformal \cite{Deshmukh}, if it satisfies $L_{\xi} g_N = 2fg_N$, where $L_{\xi}$ is the Lie derivative of the metric tensor of $g_N$ with respect to $\xi$ and $f$ is called the potential function of $\xi$.
	
	In \cite{Pigola}, Pigola introduced a natural extension of the concept of gradient Ricci soliton by taking $\lambda$ as a variable function instead of a constant and then the Ricci soliton $(M,g_1,\xi,\lambda)$ is called an almost Ricci soliton. Hence, the almost Ricci soliton becomes a Ricci soliton, if the function $\lambda$ is a constant. The almost Ricci soliton is called shrinking, steady or expanding accordingly as $\lambda<0$, $\lambda=0$ or $\lambda>0$, respectively. In \cite{Perelman}, Perelman used the Ricci soliton to solve the Poincar$\acute{e}$ conjecture, and then the geometry of Ricci solitons has been the focus of attention of many mathematicians. Moreover, Ricci solitons have been studied on contact, paracontact, almost co-K\"ahler, normal almost contact and Sasakian manifolds \cite{De, Majhi, Suh}. In \cite{Meric2019}, Meri\c{c} and Kili\c{c} studied Riemannian submersions whose total manifolds admit a Ricci soliton. In \cite{Siddiqi}, Siddiqi and Akyol studied $\eta$-Ricci-Yamabe solitons on Riemannian submersions from Riemannian manifolds. In \cite{Meric2020}, Meri\c{c} studied the Riemannian submersions admitting an almost Yamabe soliton.
	Recently, present authors introduced Riemannian maps and conformal submersions whose total manifolds admit a Ricci soliton in \cite{Yadav2021}, \cite{Yadav2022} and \cite{Meena}.
	
	In this paper, we study Riemannian maps whose base manifolds admit a Ricci soliton. In section \ref{sec2}, we recall some basic facts on Riemannian maps which are needed for this paper. In section \ref{sec3}, a Riemannian map $F$ between Riemannian manifold is considered and we find the Riemannian curvature tensor of the base manifold. Moreover, we calculate the Ricci tensor and the scalar curvature of the base manifold. In section \ref{sec4}, we obtain necessary conditions for leaves of $rangeF_\ast$ to be Ricci soliton, almost Ricci soliton and Einstein. We also obtain necessary conditions for leaves of $(rangeF_\ast)^\bot$ to be Ricci soliton and Einstein. Moreover, we calculate the scalar curvatures of $rangeF_\ast$ and $(rangeF_\ast)^\bot$ for a totally geodesic Riemannian map $F$ by using Ricci soliton. Section \ref{sec5} is devoted to harmonicity and biharmonicity, in which we obtain a necessary and sufficient condition for a Riemannian map between Riemannian manifolds whose base manifold admits a Ricci soliton to be harmonic. We also obtain necessary and sufficient conditions for a Riemannian map from a Riemannian manifold to a space form which admits Ricci soliton to be harmonic and biharmonic. In the last section, we give a non-trivial example of a Riemannian map whose base manifold admits a Ricci soliton.
	\section{Preliminaries}\label{sec2}
	In this section, we recall the notion of Riemannian maps between Riemannian manifolds and give a brief review of basic facts of Riemannian maps.
	
	Let $ F:(M^m,g_M)\rightarrow (N^n,g_N)$ be a smooth map between Riemannian manifolds such that $ 0 < rankF\leq min\{m,n\}$, where $ dim (M)=m $ and $ dim (N)=n$. Then we denote the kernel space of $F_\ast$ by $\nu_p= kerF_{\ast p}$ at $p\in M$ and consider the orthogonal complementary space $\mathcal{H}_p=(kerF_{\ast p})^\bot$ to $kerF_{\ast p}$ in $T_pM$. Then the tangent space $T_pM$ of $M$ at $p$ has the decomposition $T_pM= (kerF_{\ast p}) \oplus(kerF_{\ast p})^\bot = \nu_p \oplus \mathcal{H}_p$. We denote the range of $F_\ast $ by $rangeF_{\ast p}$ at $p \in M$ and consider the orthogonal complementary space $(rangeF_{\ast p})^\bot$ to  $rangeF_{\ast p}$ in the tangent space $T_{F(p)}N$ of $N$ at $F(p) \in N$. If $rankF < min \{m,n \}$, then we have $(rangeF_\ast)^\bot \neq \{0\}$. Thus the tangent space $T_{F(p)}N$ of $N$ at $F(p) \in N$ has the decomposition $T_{F(p)}N= (rangeF_{\ast p}) \oplus(rangeF_{\ast p})^\bot$. 
	
	Now, a smooth map $ F:(M^m,g_M) \rightarrow (N^n,g_N)$ is called Riemannian map at $p\in M$ if the horizontal restriction $F^h_{\ast p}:(kerF_{\ast p})^\bot \rightarrow (rangeF_{\ast p})$ is a linear isometry between the inner product spaces $((kerF_{\ast p})^\bot,g_{M(p)}|_{(kerF_{\ast p})^\bot} ) $ and $(rangeF_{\ast p},g_{N(p_1)}|_{(rangeF_{\ast p})})$, where $ F(p)= p_1 $. In other words, $F_\ast$ satisfies the equation
	\begin{equation}\label{eqn2.1} 
		g_N(F_\ast X, F_\ast Y) =g_M(X,Y),
	\end{equation}
	for all $X,Y$ vector fields tangent to $\Gamma(kerF_{\ast p})^\bot$. It follows that isometric immersions and Riemannian submersions are particular Riemannian maps with $kerF_\ast= \{0\}$ and $(rangeF_\ast)^\bot= \{0\}$, respectively.
	
	Let $F:(M,g_M)\rightarrow (N,g_N)$ be a smooth map between Riemannian manifolds. Then the differential $F_\ast$ of $F$ can be viewed as a section of bundle $Hom(TM,F^{-1}TN)$ $\rightarrow M$, where $F^{-1}TN$ is the pullback bundle whose fibers at $p\in M$ is $(F^{-1}TN)_p = T_{F(p)}N$, $p\in M$. The bundle $Hom(TM,F^{-1}TN)$ has a connection $\nabla$ induced from the Levi-Civita connection ${\nabla}^M$ and the pullback connection  $\overset{N}{\nabla^F}$. Then the second fundamental form of $F$ is given by
	\begin{equation}\label{eqn2.2} 
		(\nabla F_\ast) (X,Y) = \overset{N}{\nabla_X^F} F_\ast Y - F_\ast({\nabla}_X^M Y),
	\end{equation}
	for all $X,Y \in \Gamma(TM)$, where $\overset{N}{\nabla_{X}^F} F_\ast Y \circ F= \nabla_{F_\ast X}^N F_\ast Y$. It is known that the second fundamental form is symmetric. In \cite{Sahin2010}, B. \c{S}ahin proved that $(\nabla F_\ast)(X,Y)$ has no component in $rangeF_\ast$ for all $ X,Y \in \Gamma(kerF_\ast)^\bot$. More precisely, we have
	\begin{equation}\label{eqn2.3} 
		(\nabla F_\ast) (X,Y)\in \Gamma(rangeF_\ast)^\bot.
	\end{equation}
	For any vector field $X$ on $M$ and any section $V$ of $(rangeF_\ast)^\bot,$ we have $\nabla_X^{F\bot} V$, which is the orthogonal projection of $\nabla_X^N V$ on $(rangeF_\ast)^\bot,$ where $\nabla^{F \bot}$ is a linear connection on $(rangeF_\ast)^\bot$ such that $\nabla^{F \bot} g_N = 0$.
	
	\noindent Now, for a Riemannian map $F$ we define $\mathcal{S}_V$ as (\cite{Sahinbook}, p. 188)
	\begin{equation}\label{eqn2.4} 
		\nabla_{F_\ast X}^N V = -\mathcal{S}_V F_\ast X + \nabla_X^{F \bot} V,
	\end{equation}
	where $\nabla^N$ is the Levi-Civita connection on $N$, $\mathcal{S}_V F_\ast X$ is the tangential component (a vector field along $F$) of $\nabla_{F_\ast X}^N V.$ Thus at $p\in M$, we have $\nabla_{F_\ast X}^N V(p) \in T_{F(p)} N$, $\mathcal{S}_V F_\ast X \in F_{\ast p} (T_p M)$ and $\nabla_{X}^{F \bot} V(p) \in (F_{\ast p}(T_p M))^\bot.$ It is easy to see that $\mathcal{S}_V F_\ast X$ is bilinear in $V$, and $F_\ast X$ at $p$ depends only on $V_p$ and $F_{\ast p} X_p.$ Hence from (\ref{eqn2.2}) and (\ref{eqn2.4}), we obtain
	\begin{equation}\label{eqn2.5} 
		g_N(\mathcal{S}_V F_\ast X, F_\ast Y) = g_N(V, (\nabla F_\ast)(X,Y)),
	\end{equation}
	for $X,Y \in \Gamma(kerF_\ast)^\bot$ and $V \in \Gamma(rangeF_\ast)^\bot.$ Using (\ref{eqn2.5}), we obtain
	\begin{equation}\label{eqn2.6} 
		g_N((\nabla F_\ast)(X, {}^\ast F_\ast(\mathcal{S}_V F_\ast Y)), W) = g_N(\mathcal{S}_W F_\ast X, \mathcal{S}_V F_\ast Y),
	\end{equation}
	where $\mathcal{S}_V$ is self-adjoint operator and ${}^\ast F_\ast$ is the adjoint map of $F_\ast$. For details, we refer to (\cite{Sahinbook}, p. 186).
	
	\begin{definition}
		\cite{Sahin2011} A Riemannian map $F$ between Riemannian manifolds $(M,g_M)$ and $(N,g_N)$ is said to be an umbilical Riemannian map at $p \in M$, if
		\begin{equation}\label{eqn2.7} 
			\mathcal{S}_V F_{\ast, p}(X_p) = f F_{\ast p} X_p,
		\end{equation}
		for any $F_\ast X \in \Gamma(rangeF_\ast)$ and $V \in \Gamma(rangeF_\ast)^\bot,$ where $f$ is a differential function on $M$. If $F$ is umbilical for every $p \in M$ then we say that $F$ is an umbilical Riemannian map.
	\end{definition} 
   
	The Riemannian curvature tensor $R^N$ of $N$ is a (1, 3) tensor field defined by $R^N (X_1, Y_1)Z_1 = \nabla_{X_1} \nabla_{Y_1} Z_1 - \nabla_{Y_1} \nabla_{X_1} Z_1 - \nabla_{[X_1,Y_1]} Z_1$ for any $X_1,Y_1,Z_1 \in \Gamma(TN)$. Now for $F_\ast X, F_\ast Y, F_\ast Z \in \Gamma(rangeF_\ast)$ and $V \in \Gamma(rangeF_\ast)^\bot,$ we have \cite{Sahin2015}
	\begin{equation}\label{eqn2.8} 
		\begin{array}{ll}
			g_N(R^N (F_\ast X, F_\ast Y)V, F_\ast Z) = &g_N((\tilde{\nabla}_Y \mathcal{S})_V F_\ast X, F_\ast Z) \\&- g_N((\tilde{\nabla}_X \mathcal{S})_V F_\ast Y, F_\ast Z),
		\end{array}
	\end{equation}
	where $(\tilde{\nabla}_X \mathcal{S})_V F_\ast Y$ is defined by
	\begin{equation}\label{eqn2.9} 
		(\tilde{\nabla}_X \mathcal{S})_V F_\ast Y = F_\ast ( \nabla_X^M {}^\ast F_\ast (\mathcal{S}_V F_\ast Y)) -  \mathcal{S}_{\nabla_{X}^{F \bot} V} F_\ast Y - \mathcal{S}_V P \overset{N}{\nabla_X^{F}} F_\ast Y,
	\end{equation}
	where $P$ denotes the projection morphism on $rangeF_\ast$ and $R^N$ is the Riemannian curvature tensor of $\nabla^N$ (which is a metric connection on $N$).
	
	If $N$ is of constant sectional curvature $c$, denoted by $N(c)$ (known as space form), whose curvature tensor field $R^N$ is given by \cite{Yano}
	\begin{equation}\label{eqn2.10} 
		R^N(X_1,Y_1)Z_1 = c \{g_N(Y_1,Z_1)X_1 - g_N(X_1,Z_1)Y_1\},
	\end{equation}
	for $X_1,Y_1,Z_1 \in \Gamma(TN)$. Now we denote the Ricci tensor and the scalar curvature by $Ric$ and $s^N$, respectively and defined as $Ric (X_1,Y_1)= trace (Z_1 \mapsto R(Z_1,X_1)Y_1)$ and $s^N= trace Ric (X_1,Y_1)$ for $ X_1,Y_1 \in \Gamma(TN)$.
	
	The gradient of a smooth function $f$ denoted by $gradf$ and is defined as 
	\begin{equation}\label{eqn2.11} 
		g_N (gradf, X_1) = X_1(f),
	\end{equation}
	for $X_1 \in \Gamma(TN)$.
	\section{The equations of Riemannian curvature, Ricci tensor and scalar curvature for the base manifold of a Riemannian map}\label{sec3}
	In this section, we will find Riemannian curvature tensor and then calculate the Ricci tensor and the scalar curvature of base manifold.
	 
	 Given a Riemannian map $F:(M,g_M) \rightarrow (N,g_N)$ we denote ${}^\ast F_\ast$ the adjoint map of $F_\ast$. For any $X \in \Gamma(TM), V \in \Gamma(rangeF_\ast)^\bot$ we denote by $\nabla_X^{F \bot} V$ the projection of $\nabla_X^N V$ on $(rangeF_\ast)^\bot$. This allows to define a linear connection $\nabla^{F \bot}$ on $(range F_\ast)^\bot$. Now, for all $U, V \in \Gamma(rangeF_\ast)^\bot$, we define
	\begin{equation*}
		\nabla^N_U V = \mathcal{R}(\nabla^N_U V) + \nabla^{F \bot}_U V,
	\end{equation*}
	where $\mathcal{R}(\nabla^N_U V)$ and $\nabla^{F \bot}_U V$ denote the component of $\nabla^N_U V$ on $rangeF_\ast$ and $(rangeF_\ast)^\bot$, respectively. Then the distribution $(rangeF_\ast)^\bot$ is totally geodesic if and only if $\mathcal{R}(\nabla^N_U V) = 0$.
	Note that throughout this paper, we assumed the Riemannian map $F:(M,g_M) \rightarrow (N,g_N)$ such that {\small $(rangeF_\ast)^\bot$} is totally geodesic, that is $\nabla^N_U V = \nabla^{F \bot}_U V~ \text{for all}~ U, V \in \Gamma(rangeF_\ast)^\bot$.
	
	\begin{proposition}  Let $F:(M,g_M) \rightarrow (N,g_N)$ be a Riemannian map between Riemannian manifolds. Then for any $X\in \Gamma(TM)$, $V,W \in \Gamma(rangeF_\ast)^\bot$, we have
		\begin{equation}\label{eqn3.1} 
			\begin{array}{ll}
				R^N(F_\ast X, V)W=& - \mathcal{S}_{\nabla_V^{F \bot} W} F_\ast X + \nabla_X^{F \bot} \nabla_V^{F \bot} W + \nabla_V^N \mathcal{S}_W F_\ast X \\& - \nabla_V^{F \bot} \nabla_X^{F \bot} W - \mathcal{S}_W\mathcal{S}_V F_\ast X + \nabla_{ {}^\ast F_\ast(\mathcal{S}_V F_\ast X)}^{F \bot} W \\&- \nabla_{\nabla_X^{F \bot} V}^{F \bot} W - \mathcal{S}_W \nabla_V^N F_\ast X + \nabla_{ {}^\ast F_\ast(\nabla_V^N F_\ast X)}^{F \bot} W.
			\end{array}
		\end{equation}
		The component of $R^N(F_\ast X, V)W$ on $(rangeF_\ast)^\bot$  is
			\begin{equation}\label{eqn3.2} 
			\begin{array}{ll}
				R^{F \bot}(F_\ast X, V)W=& \nabla_X^{F \bot} \nabla_V^{F \bot} W - \nabla_V^{F \bot} \nabla_X^{F \bot} W- \nabla_{\nabla_X^{F \bot} V}^{F \bot} W  \\&+ \nabla_{ {}^\ast F_\ast(\mathcal{S}_V F_\ast X)}^{F \bot} W + \nabla_{ {}^\ast F_\ast(\nabla_V^N F_\ast X)}^{F \bot} W.
			\end{array}
		\end{equation}
	\end{proposition}
	\begin{proof} Let $ F:(M,g_M)\rightarrow (N,g_N)$ be a Riemannian map between Riemannian manifolds. Now for $F_\ast X \in \Gamma(rangeF_\ast)$ and $V,W \in \Gamma(rangeF_\ast)^\bot,$ we have  
		\begin{equation}\label{eqn3.3} 
			R^N(F_\ast X, V)W= \nabla_{F_\ast X}^N \nabla_V^N W - \nabla_V^N \nabla_{F_\ast X}^N W - \nabla_{[F_\ast X, V]}^{N} W.
		\end{equation}
		Now, using (\ref{eqn2.4}), we get
		\begin{equation}\label{eqn3.4} 
			\nabla_{F_\ast X}^N \nabla_V^N W = \nabla_{F_\ast X}^N \nabla_V^{F \bot} W = - \mathcal{S}_{\nabla_V^{F \bot} W} F_\ast X + \nabla_X^{F \bot} \nabla_V^{F \bot} W,
		\end{equation}
		\begin{equation}\label{eqn3.5} 
			\nabla_V^N \nabla_{F_\ast X}^N W = -\nabla_V^N \mathcal{S}_W F_\ast X + \nabla_V^{F \bot} \nabla_X^{F \bot} W,
		\end{equation}
		and
		\begin{equation}\label{eqn3.6} 
			\nabla_{[F_\ast X, V]}^{N} W= \mathcal{S}_W\mathcal{S}_V F_\ast X - \nabla_{ {}^\ast F_\ast(\mathcal{S}_V F_\ast X)}^{F \bot} W + \nabla_{\nabla_X^{F \bot} V}^{F \bot} W - \nabla_{\nabla_V^N F_\ast X}^N W.
		\end{equation}
		Since $g_N(\nabla_V^N F_\ast X, U)$ $= 0$ for all $U\in \Gamma(rangeF_\ast)^\bot$, $\nabla_V^N F_\ast X \in \Gamma(rangeF_\ast)$. Then using (\ref{eqn2.4}) in (\ref{eqn3.6}), we get
		\begin{equation}\label{eqn3.7} 
			\begin{array}{ll}
				\nabla_{[F_\ast X, V]}^{N} W=& \mathcal{S}_W\mathcal{S}_V F_\ast X - \nabla_{ {}^\ast F_\ast(\mathcal{S}_V F_\ast X)}^{F \bot} W + \nabla_{\nabla_X^{F \bot} V}^{F \bot} W \\&+ \mathcal{S}_W \nabla_V^N F_\ast X - \nabla_{ {}^\ast F_\ast(\nabla_V^N F_\ast X)}^{F \bot} W.
			\end{array}
		\end{equation}
		Now using (\ref{eqn3.4}), (\ref{eqn3.5}) and (\ref{eqn3.7}) in (\ref{eqn3.3}), we get (\ref{eqn3.1}). This completes the proof.
	\end{proof}
	
	\noindent Now, we examine the following consequences of Proposition 3.1.
	\begin{lemma} 
		 Let $F:(M,g_M) \rightarrow (N,g_N)$ be a Riemannian map between \\Riemannian manifolds. Then for any $F_\ast X, F_\ast Y \in \Gamma(rangeF_\ast)$ and $U, V, W \in  \Gamma(rangeF_\ast)^\bot$, we have
	\begin{equation}\label{eqn3.8} 	
		\begin{array}{ll}
			g_N( R^N(F_\ast X, V)W, F_\ast Y)=& -g_N( \mathcal{S}_{\nabla_V^{F \bot} W} F_\ast X, F_\ast Y) \\&+ g_N( \nabla_V^N \mathcal{S}_W F_\ast X, F_\ast Y) \\& -g_N( \mathcal{S}_V F_\ast X, \mathcal{S}_W F_\ast Y) \\&- g_N( \mathcal{S}_W (\nabla_V^N F_\ast X), F_\ast Y),
		\end{array}
	\end{equation}
	and
	\begin{equation}\label{eqn3.9} 
		\begin{array}{ll}
			g_N( R^N(F_\ast X, V)W, U)=& g_N\Big ( \nabla_X^{F \bot} \nabla_V^{F \bot} W- \nabla_V^{F \bot} \nabla_X^{F \bot} W \\&+ \nabla_{ {}^\ast F_\ast(\mathcal{S}_V F_\ast X)}^{F \bot} W - \nabla_{\nabla_X^{F \bot} V}^{F \bot} W \\&+ \nabla_{ {}^\ast F_\ast(\nabla_V^N F_\ast X)}^{F \bot} W, U\Big).
		\end{array}
	\end{equation}
	\end{lemma}
	\begin{theorem} Let $ F:(M^m,g_M)\rightarrow (N^n,g_N)$ be a Riemannian map between Riemannian manifolds. Then, the Ricci tensor on $(N,g_N)$ acts as
		\begin{equation}\label{eqn3.10} 
			\begin{array}{ll}
				Ric(F_\ast X, F_\ast Y)=&  Ric^{rangeF_\ast}(F_\ast X, F_\ast Y) - \sum\limits_{k=1}^{n_1} \Big \{ g_N(\mathcal{S}_{\nabla_{e_k}^{F \bot} e_k} F_\ast X, F_\ast Y) \\&-  g_N(\nabla_{e_k}^N \mathcal{S}_{e_k} F_\ast X, F_\ast Y)+g_N(\mathcal{S}_{e_k}  F_\ast X, \mathcal{S}_{e_k} F_\ast Y)  \\&+  g_N(\nabla_{e_k}^N  F_\ast X, \mathcal{S}_{e_k} F_\ast Y) \Big \},
			\end{array}
		\end{equation}
		\begin{equation}\label{eqn3.11} 
			\begin{array}{ll}
				Ric(V, W)=&  Ric^{(rangeF_\ast)^\bot}(V,W) - \sum\limits_{j=r+1}^{m}\Big \{ g_N(\mathcal{S}_{\nabla_{V}^{F \bot} W} F_\ast X_j, F_\ast X_j) \\&+  g_N(\mathcal{S}_{V} F_\ast X_j, \mathcal{S}_{W} F_\ast X_j)- \nabla_V^N (g_N(\mathcal{S}_{W}  F_\ast X_j, F_\ast X_j))  \\&+  2 g_N( \mathcal{S}_{W} F_\ast X_j, \nabla_{V}^N  F_\ast X_j) \Big \},
			\end{array}
		\end{equation}
		and
		\begin{equation}\label{eqn3.12} 
			\begin{array}{ll}
				Ric(F_\ast X,V)=& \sum\limits_{j=r+1}^{m}\Big \{ g_N((\tilde{\nabla}_X \mathcal{S})_V F_\ast X_j, F_\ast X_j) \\&- g_N((\tilde{\nabla}_{X_j} \mathcal{S})_V F_\ast X, F_\ast X_j)\Big \} \\&-  \sum\limits_{k=1}^{n_1} g_N \Big ( \nabla_X^{F \bot} \nabla_{e_k}^{F \bot} V- \nabla_{e_k}^{F \bot}\nabla_X^{F \bot} V - \nabla_{\nabla_X^{F \bot} e_k}^{F \bot} V \\&+ \nabla_{ {}^\ast F_\ast(\mathcal{S}_{e_k} F_\ast X)}^{F \bot} V + \nabla_{ {}^\ast F_\ast(\nabla_{e_k}^N F_\ast X)}^{F \bot} V, e_k \Big ),
			\end{array}
		\end{equation}
		for $X,Y \in \Gamma(kerF_\ast)^\bot$, $V,W \in \Gamma(rangeF_\ast)^\bot$ and $F_\ast X, F_\ast Y \in \Gamma(rangeF_\ast)$, where $\{F_\ast X_j\}_{r+1 \leq j \leq m}$ and $\{ e_k\}_{1\leq k \leq n_1}$ are orthonormal bases of $rangeF_\ast$ and $(rangeF_\ast)^\bot$, respectively and ${}^\ast F_\ast$ is the adjoint map of $F_\ast$.
	\end{theorem}
	\begin{proof}  We know that
		\begin{equation*}
			\begin{array}{ll}
				Ric(F_\ast X, F_\ast Y)=& \sum\limits_{j=r+1}^{m}g_N(R^N(F_\ast X_j, F_\ast X) F_\ast Y, F_\ast X_j) \\&+ \sum\limits_{k=1}^{n_1} g_N(R^N (e_k, F_\ast X) F_\ast Y, e_k),
			\end{array}
		\end{equation*}
		for $X,Y \in \Gamma(kerF_\ast)^\bot$, where $\{F_\ast X_j\}_{r+1 \leq j \leq m}$ and $\{ e_k\}_{1\leq k \leq n_1}$ are orthonormal bases of $rangeF_\ast$ and $(rangeF_\ast)^\bot$, respectively. Then using (\ref{eqn3.8}) in the above equation, we get (\ref{eqn3.10}).
		
		\noindent Also, we know that
		\begin{equation*}
			Ric(V, W)= \sum\limits_{j=r+1}^{m}g_N(R^N(F_\ast X_j, V) W, F_\ast X_j) + \sum\limits_{k=1}^{n_1} g_N(R^N (e_k, V) W, e_k),
		\end{equation*}
		for $V,W \in \Gamma(rangeF_\ast)^\bot$. Then using (\ref{eqn3.8}) in above equation, we get
		\begin{equation}\label{eqn3.13} 
			\begin{array}{ll}
				Ric(V, W)=&  Ric^{(rangeF_\ast)^\bot}(V,W) + \sum\limits_{j=r+1}^{m}\Big \{ - g_N(\mathcal{S}_{\nabla_{V}^{F \bot} W} F_\ast X_j, F_\ast X_j) \\&+ g_N(\nabla_{V}^N \mathcal{S}_{W} F_\ast X_j,   F_\ast X_j)  -  g_N(\mathcal{S}_{V} F_\ast X_j, \mathcal{S}_{W} F_\ast X_j) \\&- g_N(\mathcal{S}_{W}\nabla_V^N   F_\ast X_j, F_\ast X_j)\Big \}.
			\end{array}
		\end{equation}
		Since $\nabla^N$ is a metric connection on $N$, by (\ref{eqn3.13}) we get (\ref{eqn3.11}). Similarly, by using (\ref{eqn2.8}) and (\ref{eqn3.9}), we get (\ref{eqn3.12}). This completes the proof.
	\end{proof}
	\begin{theorem}\label{thm3.2} Let $ F:(M^m,g_M)\rightarrow (N^n,g_N)$ be a Riemannian map between Riemannian manifolds. Then
		\begin{equation*}
			\begin{array}{ll}
				s^N=& s^{rangeF_\ast} + s^{(rangeF_\ast)^\bot} \\&-2 \sum\limits_{j=r+1}^{m}\sum\limits_{k=1}^{n_1} g_N(\mathcal{S}_{\nabla_{e_k}^{F \bot} e_k} F_\ast X_j, F_\ast X_j) \\&+ \sum\limits_{j=r+1}^{m}\sum\limits_{k=1}^{n_1}g_N(\nabla_{e_k}^N \mathcal{S}_{e_k}F_\ast X_j,  F_\ast X_j)\\&-2\sum\limits_{j=r+1}^{m}\sum\limits_{k=1}^{n_1} g_N(\mathcal{S}_{e_k} F_\ast X_j, \mathcal{S}_{e_k} F_\ast X_j ) \\&- 3 \sum\limits_{j=r+1}^{m}\sum\limits_{k=1}^{n_1}g_N(\nabla_{e_k}^N F_\ast X_j,\mathcal{S}_{e_k}  F_\ast X_j) \\&+\sum\limits_{j=r+1}^{m}\sum\limits_{k=1}^{n_1} \nabla_{e_k}^N(g_N( \mathcal{S}_{e_k}F_\ast X_j,  F_\ast X_j)) ,
			\end{array}
		\end{equation*}
		where $s^N, s^{rangeF_\ast}$ and $s^{(rangeF_\ast)^\bot}$ denote the scalar curvatures of $N$, $rangeF_\ast$ and $(rangeF_\ast)^\bot$, respectively. In addition $\{F_\ast X_j\}_{r+1 \leq j \leq m}$ and $\{e_k\}_{1 \leq k \leq n_1}$ are orthonormal bases of $rangeF_\ast$ and $(rangeF_\ast)^\bot$.
	\end{theorem}
	\begin{proof}  Since scalar curvature of $N$ is defined by
		\begin{equation}\label{eqn3.14} 
			s^N	= \sum_{l=r+1}^{m} Ric(F_\ast X_l, F_\ast X_l) + \sum_{t=1}^{n_1} Ric(e_t, e_t),
		\end{equation}
		where $\{F_\ast X_l\}_{r+1 \leq l \leq m}$ and $\{e_t\}_{1 \leq t \leq n_1}$ are orthonormal bases of $rangeF_\ast$ and $(rangeF_\ast)^\bot$, respectively.
		Now, using (\ref{eqn3.10}) and (\ref{eqn3.11}) in (\ref{eqn3.14}), we get
		\begin{equation*}
			\begin{array}{ll}
				s^N=& \sum\limits_{l=r+1}^{m}\sum\limits_{k=1}^{n_1} \Big \{ Ric^{rangeF_\ast}(F_\ast X_l, F_\ast X_l)- g_N(\mathcal{S}_{\nabla_{e_k}^{F \bot} e_k} F_\ast X_l, F_\ast X_l) \\&+ g_N(\nabla_{e_k}^N \mathcal{S}_{e_k}F_\ast X_l,  F_\ast X_l)- g_N(\mathcal{S}_{e_k} F_\ast X_l, \mathcal{S}_{e_k} F_\ast X_l ) \\&-g_N(\nabla_{e_k}^N F_\ast X_l,\mathcal{S}_{e_k}  F_\ast X_l) \Big \} + \sum\limits_{j=r+1}^{m}\sum\limits_{t=1}^{n_1} \Big \{ Ric^{(rangeF_\ast)^\bot}(e_t, e_t) \\&- g_N(\mathcal{S}_{\nabla_{e_t}^{F \bot} e_t} F_\ast X_j, F_\ast X_j) - g_N(\mathcal{S}_{e_t} F_\ast X_j, \mathcal{S}_{e_t} F_\ast X_j ) \\&+ \nabla_{e_t}^N(g_N( \mathcal{S}_{e_t}F_\ast X_j,  F_\ast X_j)) -2 g_N(\mathcal{S}_{e_t}  F_\ast X_j, \nabla_{e_t}^N F_\ast X_j) \Big \},
			\end{array}
		\end{equation*}
		which implies
		\begin{equation*}
			\begin{array}{ll}
				s^N=& s^{rangeF_\ast} + s^{(rangeF_\ast)^\bot} \\&-2 \sum\limits_{j=r+1}^{m}\sum\limits_{k=1}^{n_1} g_N(\mathcal{S}_{\nabla_{e_k}^{F \bot} e_k} F_\ast X_j, F_\ast X_j) \\&+ \sum\limits_{j=r+1}^{m}\sum\limits_{k=1}^{n_1}g_N(\nabla_{e_k}^N \mathcal{S}_{e_k}F_\ast X_j,  F_\ast X_j)\\&-2\sum\limits_{j=r+1}^{m}\sum\limits_{k=1}^{n_1} g_N(\mathcal{S}_{e_k} F_\ast X_j, \mathcal{S}_{e_k} F_\ast X_j ) \\&- 3 \sum\limits_{j=r+1}^{m}\sum\limits_{k=1}^{n_1}g_N(\nabla_{e_k}^N F_\ast X_j,\mathcal{S}_{e_k}  F_\ast X_j) \\&+\sum\limits_{j=r+1}^{m}\sum\limits_{k=1}^{n_1} \nabla_{e_k}^N(g_N( \mathcal{S}_{e_k}F_\ast X_j,  F_\ast X_j)).
			\end{array}
		\end{equation*}
		This completes the proof.
	\end{proof}
	\begin{corollary} Let $ F:(M^m,g_M)\rightarrow (N^n,g_N)$ be a totally geodesic Riemannian map between Riemannian manifolds. Then
		\begin{equation*}
			s^N= s^{rangeF_\ast} + s^{(rangeF_\ast)^\bot}.
		\end{equation*}
	\end{corollary}
	\begin{proof} Since $F$ is totally geodesic then $\mathcal{S}_V F_\ast X = 0$ for all $X\in \Gamma(kerF_\ast)^\bot$ and $V\in \Gamma(rangeF_\ast)^\bot$. Then, the statement follows by Theorem \ref{thm3.2}.
	\end{proof}
	\begin{corollary} Let $ F:(M^m,g_M)\rightarrow (N^n,g_N)$ be an umbilical Riemannian map between Riemannian manifolds. Then
		\begin{equation*}
			\begin{array}{ll}
				s^N= s^{rangeF_\ast} + s^{(rangeF_\ast)^\bot} -2 ( f+f^2 )(m-r).
			\end{array}
		\end{equation*}
	\end{corollary}
	\begin{proof} Since $F$ is an umbilical map then using (\ref{eqn2.7}) in Theorem \ref{thm3.2}, we get
		\begin{equation*}
			\begin{array}{ll}
				s^N=& s^{rangeF_\ast} + s^{(rangeF_\ast)^\bot} -2 \sum\limits_{j=r+1}^{m}\sum\limits_{k=1}^{n_1} g_N(f F_\ast X_j, F_\ast X_j) \\&+ \sum\limits_{j=r+1}^{m}\sum\limits_{k=1}^{n_1}g_N(\nabla_{e_k}^N f F_\ast X_j,  F_\ast X_j)-2\sum\limits_{j=r+1}^{m}\sum\limits_{k=1}^{n_1} g_N(f F_\ast X_j, f F_\ast X_j ) \\&-3  \sum\limits_{j=r+1}^{m}\sum\limits_{k=1}^{n_1}g_N(\nabla_{e_k}^N F_\ast X_j,f  F_\ast X_j) +\sum\limits_{j=r+1}^{m}\sum\limits_{k=1}^{n_1} \nabla_{e_k}^N(g_N( f F_\ast X_j,  F_\ast X_j)),
			\end{array}
		\end{equation*}
		which implies the proof.
	\end{proof}
	\section{Riemannian map whose base manifold admits a Ricci soliton}\label{sec4}
	In this section, we consider a Riemannian map $ F:(M,g_M) \rightarrow (N,g_N)$ from a Riemannian manifold to a Ricci soliton and give some characterizations.
	\begin{proposition}\label{prop4.1}  \cite{Sahin2012}  Let $F:(M,g_M) \rightarrow (N,g_N)$ be a Riemannian map between Riemannian manifolds. Then $F$ is totally geodesic if and only if\\ $(i)$ $A_X Y = 0$,\\ $(ii)$ the fibers of $F$ define totally geodesic foliation on $M$,\\ $(iii)$ $S_V F_\ast X = 0$,\\for $X,Y \in \Gamma(kerF_\ast)^\bot$ and $V \in \Gamma(rangeF_\ast)^\bot$.
	\end{proposition}
	\begin{remark}
		 Since $rangeF_\ast$ is a subbundle of $TN$, it defines a distribution on $N$. Then for $F_\ast X, F_\ast Y \in \Gamma(rangeF_\ast)$, we have
		\begin{align*}
			[F_\ast X, F_\ast Y]& = \nabla_{F_\ast X}^N F_\ast Y - \nabla_{F_\ast Y}^N F_\ast X\\& = \overset{N}{\nabla_X^F} F_\ast Y \circ F -\overset{N}{\nabla_Y^F} F_\ast X \circ F.
		\end{align*}
		Using (\ref{eqn2.2}) in above equation, we get
		\begin{align*}
			[F_\ast X, F_\ast Y]& =  F_\ast(\nabla_X Y) - F_\ast(\nabla_Y X) = F_\ast(\nabla_X Y - \nabla_Y X) \in \Gamma(rangeF_\ast).
		\end{align*}
	Thus $rangeF_\ast$ is an integrable distribution. Then for any point $F(p) \in N$ there exists maximal integral manifold or a leaf of $rangeF_\ast$ containing $F(p)$.
	\end{remark}

	\begin{theorem} Let $F:(M,g_M) \rightarrow (N,g_N)$ be a totally geodesic Riemannian map between Riemannian manifolds and $(N,g_N,\xi,\lambda)$ be a Ricci soliton with potential vector field $\xi \in \Gamma(TN)$. Then the following statements are true:\\ $(i)$ If the vector field $\xi = F_\ast Z$(say) $\in \Gamma(rangeF_\ast)$ with $Z \in \Gamma(kerF_\ast)^\bot$, then any leaf of $rangeF_\ast$ is a Ricci soliton.\\$(ii)$ If the vector field $\xi = V$(say) $\in \Gamma(rangeF_\ast)^\bot$, then any leaf of $rangeF_\ast$ is an Einstein.
	\end{theorem}
	\begin{proof} Since $(N,g_N,\xi,\lambda)$ is a Ricci soliton then, we have
		\begin{equation}\label{eqn4.1} 
			\frac{1}{2} (L_{\xi} g_N)(F_\ast X,F_\ast Y) + Ric (F_\ast X, F_\ast Y) + \lambda g_N(F_\ast X, F_\ast Y) = 0,
		\end{equation}
		for $F_\ast X, F_\ast Y \in \Gamma(rangeF_\ast)$. Then from (\ref{eqn4.1}), we get
		\begin{equation}\label{eqn4.2} 
			\begin{array}{ll}
				\frac{1}{2} \{g_N(\nabla_{F_\ast X}^N \xi, F_\ast Y) + g_N(\nabla_{F_\ast Y}^N \xi, F_\ast X) \} \\+ Ric (F_\ast X, F_\ast Y) + \lambda g_N(F_\ast X, F_\ast Y) = 0.
			\end{array}
		\end{equation}
		Since $F$ is totally geodesic then using $(iii)$ of Proposition \ref{prop4.1} and (\ref{eqn3.10}) in (\ref{eqn4.2}), we get
		\begin{equation}\label{eqn4.3} 
			\begin{array}{ll}
				\frac{1}{2} \{g_N(\nabla_{F_\ast X}^N \xi, F_\ast Y) + g_N(\nabla_{F_\ast Y}^N \xi, F_\ast X) \} \\+ Ric^{rangeF_\ast} (F_\ast X, F_\ast Y) + \lambda g_N(F_\ast X, F_\ast Y) = 0.
			\end{array}
		\end{equation}
		Now, if the vector field $\xi = F_\ast Z$(say) $\in \Gamma(rangeF_\ast)$, then from (\ref{eqn4.3}), we get
		\begin{equation*}
			\begin{array}{ll}
				\frac{1}{2} \{g_N(\nabla_{F_\ast X}^N F_\ast Z, F_\ast Y) + g_N(\nabla_{F_\ast Y}^N F_\ast Z, F_\ast X) \} \\+ Ric^{rangeF_\ast} (F_\ast X, F_\ast Y) + \lambda g_N(F_\ast X, F_\ast Y) = 0,
			\end{array}
		\end{equation*}
		which implies $(i)$.\\Also, if the vector field $\xi = V$(say) $\in \Gamma(rangeF_\ast)^\bot$, then from (\ref{eqn4.3}), we get
		\begin{equation*}
			\begin{array}{ll}
				\frac{1}{2} \{g_N(\nabla_{F_\ast X}^N V, F_\ast Y) + g_N(\nabla_{F_\ast Y}^N V, F_\ast X) \} \\+ Ric^{rangeF_\ast} (F_\ast X, F_\ast Y) + \lambda g_N(F_\ast X, F_\ast Y) = 0.
			\end{array}
		\end{equation*}
		Using (\ref{eqn2.4}) in above equation, we get
		\begin{equation*}
			\begin{array}{ll}
				- \frac{1}{2} \{ g_N(\mathcal{S}_V F_\ast X, F_\ast Y) + g_N(\mathcal{S}_V F_\ast Y, F_\ast X) \} \\+ Ric^{rangeF_\ast} (F_\ast X, F_\ast Y) + \lambda g_N(F_\ast X, F_\ast Y) = 0.
			\end{array}
		\end{equation*}
		Since $\mathcal{S}_V$ is self-adjoint then from above equation, we get
		\begin{equation*}
			- g_N(\mathcal{S}_V F_\ast X, F_\ast Y) + Ric^{rangeF_\ast} (F_\ast X, F_\ast Y) + \lambda g_N(F_\ast X, F_\ast Y) = 0.
		\end{equation*}
		Since $F$ is totally geodesic, again using $(iii)$ of Proposition \ref{prop4.1} in above equation, we get
		\begin{equation}\label{eqn4.4} 
			Ric^{rangeF_\ast} (F_\ast X, F_\ast Y) + \lambda g_N(F_\ast X, F_\ast Y) = 0,
		\end{equation}
		which implies $(ii)$. This completes the proof.
	\end{proof}
	\begin{theorem} Let $F:(M^m,g_M) \rightarrow (N^n,g_N)$ be a totally geodesic Riemannian map between Riemannian manifolds and $(N,g_N,\xi,\lambda)$ be a Ricci soliton with the potential vector field $\xi \in \Gamma(rangeF_\ast)^\bot$ then the scalar curvature of $rangeF_\ast$ is $-\lambda(m-r)$, where $dim(rangeF_\ast) = m-r$.
	\end{theorem}
	\begin{proof} The proof follows by (\ref{eqn4.4}).
	\end{proof}
	\begin{remark}
		Since $(rangeF_\ast)^\bot$ is a subbundle of $TN$, it defines a distribution on $N$. If $(rangeF_\ast)^\bot$ is totally geodesic then for $U, V \in \Gamma(rangeF_\ast)^\bot$, we have
		\begin{align*}
			[U, V]& = \nabla_{U}^N V - \nabla_{V}^N U\\& =\nabla_{U}^{F \bot} V - \nabla_{V}^{F \bot} U \in \Gamma(rangeF_\ast)^\bot.
		\end{align*}
		Thus $(rangeF_\ast)^\bot$ is an integrable distribution. Then for any point $F(p) \in N$ there exists maximal integral manifold or a leaf of $(rangeF_\ast)^\bot$ containing $F(p)$.
	\end{remark}
	\begin{theorem} 
		Let $F:(M,g_M) \rightarrow (N,g_N)$ be a totally geodesic Riemannian map between Riemannian manifolds and $(N,g_N,\xi,\lambda)$ be a Ricci soliton with potential vector field $\xi \in \Gamma(TN)$. Then the following statements are true:\\ $(i)$ If the vector field $\xi = V$(say) $\in \Gamma(rangeF_\ast)^\bot$, then any leaf of $(rangeF_\ast)^\bot$ is a Ricci soliton.\\$(ii)$ If the vector field $\xi = F_\ast X$(say) $\in \Gamma(rangeF_\ast)$, then any leaf of $(rangeF_\ast)^\bot$ is an Einstein.
	\end{theorem}
	\begin{proof}  Since $(N,g_N,\xi,\lambda)$ be a Ricci soliton then, we have
		\begin{equation*} 
			\frac{1}{2} (L_{\xi} g_N)(U,W) + Ric (U,W) + \lambda g_N(U,W) = 0, 
		\end{equation*}
		for $U,W \in \Gamma(rangeF_\ast)^\bot$. Then from above equation, we get
		\begin{equation*}
			\frac{1}{2} \{g_N(\nabla_U^N \xi, W) + g_N(\nabla_W^N \xi, U \} + Ric (U, W) + \lambda g_N(U, W) = 0.
		\end{equation*}
		Since $F$ is totally geodesic then using $(iii)$ of Proposition \ref{prop4.1} and (\ref{eqn3.11}) in above equation, we get
		\begin{equation}\label{eqn4.5} 
			\frac{1}{2} \{g_N(\nabla_U^N \xi, W) + g_N(\nabla_W^N \xi, U \} + Ric^{(rangeF_\ast)^\bot} (U, W) + \lambda g_N(U, W) = 0.
		\end{equation}
		Now, if the vector field $\xi = V$(say) $\in \Gamma(rangeF_\ast)^\bot$, then from (\ref{eqn4.5}), we get
		\begin{equation*}
			\frac{1}{2} \{g_N(\nabla_U^N V, W) + g_N(\nabla_W^N V, U) \} + Ric^{(rangeF_\ast)^\bot} (U, W) + \lambda g_N(U, W) = 0.
		\end{equation*}
		Since $(rangeF_\ast)^\bot$ is totally geodesic, above equation can be written as
		\begin{equation*}
			\frac{1}{2} \{g_N(\nabla_U^{F \bot} V, W) + g_N(\nabla_W^{F \bot} V, U) \} + Ric^{(rangeF_\ast)^\bot} (U, W) + \lambda g_N(U, W) = 0,
		\end{equation*}
		which implies $(i)$.\\
		Also, if the vector field $\xi = F_\ast X$(say) $\in \Gamma(rangeF_\ast)$, then from (\ref{eqn4.5}), we get
		\begin{equation*}
			\frac{1}{2} \{g_N(\nabla_U^N F_\ast X, W) + g_N(\nabla_W^N F_\ast X, U \} \\+ Ric^{(rangeF_\ast)^\bot} (U, W) + \lambda g_N(U, W) = 0.
		\end{equation*}
		Since $\nabla^N$ is metric connection, using metric compatibility condition in above equation, we get
		\begin{equation*}
			-\frac{1}{2} \{g_N(\nabla_U^N W, F_\ast X) + g_N(\nabla_W^N U, F_\ast X) \} \\+ Ric^{(rangeF_\ast)^\bot} (U, W) + \lambda g_N(U, W) = 0.
		\end{equation*}
		Since $(rangeF_\ast)^\bot$ is totally geodesic, using $\nabla_U^N W = \nabla_U^{F \bot} W$ in above equation, we get
		\begin{equation}\label{eqn4.6} 
			Ric^{(rangeF_\ast)^\bot} (U, W) + \lambda g_N(U, W) = 0,
		\end{equation}
		which implies $(ii)$. This completes the proof.
	\end{proof}
	\begin{theorem} Let $F:(M^m,g_M) \rightarrow (N^n,g_N)$ be a totally geodesic Riemannian map between Riemannian manifolds and $(N,g_N,\xi,\lambda)$ be a Ricci soliton. If the potential vector field $\xi= F_\ast X (say) \in \Gamma(rangeF_\ast)$ then the scalar curvature of $(rangeF_\ast)^\bot$ is $-\lambda n_1$, where $dim(rangeF_\ast)^\bot = n_1$.
	\end{theorem}
	\begin{proof} The proof follows by  (\ref{eqn4.6}).
	\end{proof}	
	\begin{theorem} Let $F:(M,g_M) \rightarrow (N,g_N)$ be an umbilical Riemannian map between Riemannian manifolds and $(N,g_N,\xi,\lambda)$ be a Ricci soliton with potential vector field $\xi \in \Gamma(TN)$. Then the following statements are true:\\$(i)$ If the vector field $\xi=V$(say) $\in \Gamma(rangeF_\ast)^\bot$, then any leaf of $rangeF_\ast$ is an Einstein.\\$(ii)$ If the vector field $\xi = F_\ast Z$(say) $\in \Gamma(rangeF_\ast)$, then any leaf of $rangeF_\ast$ is an almost Ricci soliton.
	\end{theorem}
	\begin{proof} Let $(N,g_N,\xi,\lambda)$ be a Ricci soliton. If $\xi=V$(say) $\in \Gamma(rangeF_\ast)^\bot$ then using (\ref{eqn2.4}) and (\ref{eqn3.10}) in (\ref{eqn4.2}), we get
		\begin{equation}\label{eqn4.7} 
			\begin{array}{ll}
				\frac{1}{2} \{ g_N(-\mathcal{S}_V F_\ast X + \nabla_X^{F \bot} V, F_\ast Y) +g_N(-\mathcal{S}_V F_\ast Y + \nabla_Y^{F \bot} V, F_\ast X)  \} \\+  Ric^{rangeF_\ast}(F_\ast X, F_\ast Y)+ \lambda g_N(F_\ast X, F_\ast Y) \\- \sum\limits_{k=1}^{n_1} \Big \{ g_N(\mathcal{S}_{\nabla_{e_k}^{F \bot} e_k} F_\ast X, F_\ast Y) -  g_N(\nabla_{e_k}^N \mathcal{S}_{e_k} F_\ast X, F_\ast Y)\\+g_N(\mathcal{S}_{e_k}  F_\ast X, \mathcal{S}_{e_k} F_\ast Y)  +  g_N(\nabla_{e_k}^N  F_\ast X, \mathcal{S}_{e_k} F_\ast Y) \Big \}  = 0.
			\end{array}
		\end{equation}
		Since $\mathcal{S}_V$ is self-adjoint then from (\ref{eqn4.7}), we get
		\begin{equation*}
			\begin{array}{ll}
				- g_N(\mathcal{S}_V F_\ast X, F_\ast Y) +  Ric^{rangeF_\ast}(F_\ast X, F_\ast Y) - \sum\limits_{k=1}^{n_1} \Big \{ g_N(\mathcal{S}_{\nabla_{e_k}^{F \bot} e_k} F_\ast X, F_\ast Y) \\-  g_N(\nabla_{e_k}^N \mathcal{S}_{e_k} F_\ast X, F_\ast Y)+g_N(\mathcal{S}_{e_k}  F_\ast X, \mathcal{S}_{e_k} F_\ast Y) +  g_N(\nabla_{e_k}^N  F_\ast X, \mathcal{S}_{e_k} F_\ast Y) \Big \} \\+ \lambda g_N(F_\ast X, F_\ast Y) = 0.
			\end{array}
		\end{equation*}
		Since $F$ is an umbilical Riemannian map then using (\ref{eqn2.7}) in above equation, we get
		\begin{equation*}
			\begin{array}{ll}
				-2f g_N( F_\ast X, F_\ast Y) +  Ric^{rangeF_\ast}(F_\ast X, F_\ast Y) \\- f^2 g_N(F_\ast X, F_\ast Y) + \lambda g_N(F_\ast X, F_\ast Y) = 0.
			\end{array}
		\end{equation*}
		Thus from above equation, we get
		\begin{equation*}
			Ric^{rangeF_\ast}(F_\ast X, F_\ast Y) -\mu g_N(F_\ast X, F_\ast Y) = 0,
		\end{equation*}
		where $\mu = 2f + f^2 - \lambda$ is a differentiable function, which implies $(i)$.
		
		\noindent Also, if $\xi=F_\ast Z$(say) $\in \Gamma(rangeF_\ast)$ then using (\ref{eqn3.10}) in (\ref{eqn4.1}), we get
		\begin{equation*}
			\begin{array}{ll}
				\frac{1}{2} (L_{F_\ast Z} g_N)(F_\ast X,F_\ast Y) +  Ric^{rangeF_\ast}(F_\ast X, F_\ast Y) \\- \sum\limits_{k=1}^{n_1} \Big \{ g_N(\mathcal{S}_{\nabla_{e_k}^{F \bot} e_k} F_\ast X, F_\ast Y) -  g_N(\nabla_{e_k}^N \mathcal{S}_{e_k} F_\ast X, F_\ast Y)\\+g_N(\mathcal{S}_{e_k}  F_\ast X, \mathcal{S}_{e_k} F_\ast Y)  +  g_N(\nabla_{e_k}^N  F_\ast X, \mathcal{S}_{e_k} F_\ast Y) \Big \} + \lambda g_N(F_\ast X, F_\ast Y) = 0.
			\end{array}
		\end{equation*}
		Since $F$ is an umbilical Riemannian map then using (\ref{eqn2.7}) in above equation, we get
		\begin{equation}\label{eqn4.8} 
			\begin{array}{ll}
				\frac{1}{2} (L_{F_\ast Z} g_N)(F_\ast X,F_\ast Y) +  Ric^{rangeF_\ast}(F_\ast X, F_\ast Y) - f g_N( F_\ast X, F_\ast Y) \\ - f^2 g_N(F_\ast X, F_\ast Y)+ \lambda g_N(F_\ast X, F_\ast Y) = 0.
			\end{array}
		\end{equation}
		Thus from (\ref{eqn4.8}), we get
		\begin{equation*}
			\begin{array}{ll}
				\frac{1}{2} (L_{F_\ast Z} g_N)(F_\ast X,F_\ast Y) +  Ric^{rangeF_\ast}(F_\ast X, F_\ast Y) \\- (f + f^2 - \lambda ) g_N(F_\ast X, F_\ast Y) = 0,
			\end{array}
		\end{equation*}
		which implies $(ii)$. This completes the proof.
	\end{proof}
	\begin{theorem} Let $F:(M,g_M) \rightarrow (N,g_N)$ be a Riemannian map from a Riemannian manifold to an Einstein manifold and $(N,g_N,\xi,\lambda)$ be a Ricci soliton with potential vector field $\xi \in \Gamma(TN)$. Then the vector field $\xi$ is killing on $N$.
	\end{theorem}
	\begin{proof} 
			Since $N$ is Einstein, $Ric (F_\ast X, F_\ast Y)$ $= -\lambda g_N(F_\ast X, F_\ast Y)$. Then by (\ref{eqn1.1}), we get
		\begin{equation*}
			\frac{1}{2} (L_{\xi} g_N)(F_\ast X,F_\ast Y)= 0,
		\end{equation*}
	for $F_\ast X, F_\ast Y \in \Gamma(rangeF_\ast)$.
	In addition, since $N$ is Einstein, $Ric (U, V)= -\lambda g_N(U, V)$. Then by (\ref{eqn1.1}), we get
		\begin{equation*}
			\frac{1}{2} (L_{\xi} g_N)(U, V)= 0,
		\end{equation*}
	for $U, V \in \Gamma(rangeF_\ast)^\bot$.	
	Similarly, we can get
	\begin{equation*}
		\frac{1}{2} (L_{\xi} g_N)(F_\ast X, V)= 0,
	\end{equation*}
	for $F_\ast X \in \Gamma(rangeF_\ast)$ and $V \in \Gamma(rangeF_\ast)^\bot$.
	This implies the proof.
	\end{proof}
	\begin{theorem} Let $F:(M,g_M) \rightarrow (N,g_N)$ be a Riemannian map between Riemannian manifolds and $(N,g_N,F_\ast U,\lambda)$ be a Ricci soliton with potential vector field $F_\ast U$ for $U\in\Gamma( kerF_\ast)$. Then $(N,g_N)$ is an Einstein manifold.
	\end{theorem}
	\begin{proof}  We know that $F_\ast U = 0$ for all $U \in \Gamma(kerF_\ast)$ and since \\$(N,g_N,F_\ast U,\lambda)$ is a Ricci soliton then from (\ref{eqn1.1}), we get 
		\begin{equation*} 
			Ric(X_1,Y_1) + \lambda g_N(X_1,Y_1) = 0,
		\end{equation*}
		for any $X_1,Y_1 \in \Gamma(TN)$, which means $N$ is an Einstein manifold. This completes the proof.
	\end{proof}
	\begin{remark}
		In \cite{Sahin2017}, B. \c{S}ahin obtained necessary and sufficient condition for the total manifold of a Riemannian map to be Einstein without using Ricci soliton. On the other hand, in above theorem, we obtain a sufficient condition for the base manifold of a Riemannian map to be Einstein using Ricci soliton.
	\end{remark}
	\section{Harmonicity and biharmonicity of Riemannian map from a Riemannian manifold to a Ricci soliton}\label{sec5}
	This section deals with the harmonicity and biharmonicity of Riemannian map from a Riemannian manifold to a Ricci soliton. 
	
	A harmonic map between Riemannian manifolds has played an important role in linking the geometry to global analysis on Riemannian manifolds as well as its importance in physics is also well established. Therefore it is an interesting question to find harmonic maps to Ricci soliton.
	We first recall that a map $F:(M^m,g_M)\rightarrow (N^n,g_N)$ between Riemannian manifolds is harmonic if and only if the tension field of $F$ vanishes at each point $p\in M$, i.e.
	\begin{equation*}
		\tau(F)= trace (\nabla F_\ast) = \sum_{i=1}^{m}(\nabla F_\ast)(e_i, e_i)= 0,
	\end{equation*}
	where $\{e_i\}_{1\leq i \leq m}$ is local orthonormal frame around a point $p \in M$ and $\nabla F_\ast$ is the second fundamental form of $F$.
	\begin{lemma}\label{lem5.1} \cite{Sahin2010a} Let $ F:(M^m,g_M)\rightarrow (N^n,g_N)$ be a Riemannian map between Riemannian manifolds. Then the tension field of $F$ is given by $\tau (F) = -rF_\ast (H) + (m-r) H_2$, where $r= dim(kerF_\ast)$, $(m-r) = rankF$, $H$ and $H_2$ are the mean curvature vector fields of the distributions $kerF_\ast$ and $rangeF_\ast$, respectively.
	\end{lemma}
	\noindent Moreover, the mean curvature vector field of $rangeF_\ast$ is defined by \cite{Sahin2010a}
	\begin{equation}\label{eqn5.1} 
		H_2 = \frac{1}{m-r}\sum_{j=r+1}^{m} \nabla_{X_j}^F F_\ast (X_j),
	\end{equation}
	where $\{X_j\}_{r+1 \leq j \leq m }$ is an orthonormal basis of $(kerF_\ast)^\bot$.
	\begin{lemma} \cite{Sahin2011} Let $ F:(M^m,g_M)\rightarrow (N^n,g_N)$ be a Riemannian map between Riemannian manifolds. Then $F$ is an umbilical map if and only if 
		\begin{equation}\label{eqn5.2} 
			(\nabla F_\ast) (X, Y) = g_M(X,Y) H_2,
		\end{equation}
		for $X,Y \in \Gamma(kerF_\ast)^\bot$ and $H_2$ is nowhere zero vector field on $(rangeF_\ast)^\bot$. 
	\end{lemma}
	\begin{theorem} Let $(N,g_N,V,\lambda)$ be a Ricci soliton with potential vector field $V \in \Gamma(rangeF_\ast)^\bot$ and $F:(M^m,g_M) \rightarrow (N^n,g_N)$ be an umbilical Riemannian map between Riemannian manifolds such that the scalar curvature of $rangeF_\ast$ is $-\lambda (m-r) \neq 0$ and $kerF_\ast$ is minimal. Then $F$ is harmonic if and only if $\sum\limits_{k=1}^{n_1} e_k=0$, where $\{ e_k\}_{1 \leq k \leq n_1}$ is an orthonormal basis of $(rangeF_\ast)^\bot$.
	\end{theorem}
	\begin{proof} By using (\ref{eqn2.6}) in (\ref{eqn4.7}), we get
		\begin{equation}\label{eqn5.3} 
			\begin{array}{ll}
				\frac{1}{2} \{ g_N(-\mathcal{S}_V F_\ast X + \nabla_X^{F \bot} V, F_\ast Y) +g_N(-\mathcal{S}_V F_\ast Y + \nabla_Y^{F \bot} V, F_\ast X)  \} \\+ Ric^{rangeF_\ast}(F_\ast X, F_\ast Y) - \sum\limits_{k=1}^{n_1} \Big \{ g_N(\mathcal{S}_{\nabla_{e_k}^{F \bot} e_k} F_\ast X, F_\ast Y) \\-  g_N(\nabla_{e_k}^N \mathcal{S}_{e_k} F_\ast X, F_\ast Y)+g_N((\nabla F_\ast)(X, {}^\ast F_\ast \mathcal{S}_{e_k} F_\ast Y), e_k) \\ +  g_N(\nabla_{e_k}^N  F_\ast X, \mathcal{S}_{e_k} F_\ast Y) \Big \} + \lambda g_N(F_\ast X, F_\ast Y) = 0,
			\end{array}
		\end{equation}
		where $\{e_k\}_{1 \leq k \leq n_1 }$ is an orthonormal basis of $(rangeF_\ast)^\bot$. Since $\mathcal{S}_V$ is self-adjoint then from (\ref{eqn5.3}), we get
		\begin{equation}\label{eqn5.4} 
			\begin{array}{ll}
				- g_N(\mathcal{S}_V F_\ast X, F_\ast Y) + Ric^{rangeF_\ast}(F_\ast X, F_\ast Y) + \lambda g_N(F_\ast X, F_\ast Y)\\- \sum\limits_{k=1}^{n_1} \Big \{ g_N(\mathcal{S}_{\nabla_{e_k}^{F \bot} e_k} F_\ast X, F_\ast Y) -  g_N(\nabla_{e_k}^N \mathcal{S}_{e_k} F_\ast X, F_\ast Y)\\+g_N((\nabla F_\ast)(X, {}^\ast F_\ast \mathcal{S}_{e_k} F_\ast Y), e_k)  +  g_N(\mathcal{S}_{e_k} \nabla_{e_k}^N  F_\ast X,  F_\ast Y) \Big \}  = 0.
			\end{array}
		\end{equation}
		Since $F$ is an umbilical Riemannian map then using (\ref{eqn2.7}) and (\ref{eqn5.2}) in (\ref{eqn5.4}), we get
		\begin{equation}\label{eqn5.5} 
			\begin{array}{ll}
				-2f g_N( F_\ast X, F_\ast Y) +  Ric^{rangeF_\ast}(F_\ast X, F_\ast Y) + \lambda g_N(F_\ast X, F_\ast Y)\\ - f \sum\limits_{k=1}^{n_1}g_N(F_\ast X, F_\ast Y) g_N(H_2, e_k)  = 0.
			\end{array}
		\end{equation}
		Taking trace of (\ref{eqn5.5}), we get
		\begin{equation*}
			\begin{array}{ll}
				-2f (m-r) + s^{rangeF_\ast}+ \lambda (m-r)  \\- f (m-r) \sum\limits_{k=1}^{n_1} g_N(H_2, e_k) = 0.
			\end{array}
		\end{equation*}
		Putting $s^{rangeF_\ast} = - \lambda (m-r)$ in above equation, we get
		\begin{equation*}
			-2f - f \sum_{k=1}^{n_1} g_N(H_2, e_k) = 0,
		\end{equation*}
		which implies
		\begin{equation*}
			\frac{-2f}{n_1} \sum\limits_{k=1}^{n_1} g_N(e_k, e_k) - f \sum\limits_{k=1}^{n_1} g_N(H_2, e_k) = 0.
		\end{equation*}
		Using (\ref{eqn2.11}) in above equation, we get
		\begin{equation*}
			\frac{-2f}{n_1} \sum\limits_{k=1}^{n_1} g_N(e_k, e_k) - f \sum\limits_{k=1}^{n_1} g_N(H_2, e_k) = 0,
		\end{equation*}
		which implies
		\begin{equation*}
			\frac{-2f}{n_1} \sum\limits_{k=1}^{n_1} e_k - f H_2 = 0.
		\end{equation*}
		Hence
		\begin{equation}\label{eqn5.6} 
			H_2 = -\frac{2}{n_1} \sum\limits_{k=1}^{n_1} e_k.
		\end{equation}
		Since $kerF_\ast$ is minimal and using (\ref{eqn5.6}) in Lemma \ref{lem5.1}, we get
		
		$\tau(F) = -(m-r) \Big \{\frac{2}{n_1} \sum\limits_{k=1}^{n_1} e_k \Big \}$, which completes the proof.
	\end{proof}
	\begin{theorem}\label{thm5.2} Let $(N,g_N,F_\ast X,\lambda)$ be a Ricci soliton with potential vector field $F_\ast X \in \Gamma(rangeF_\ast)$ and $F:(M^m,g_M) \rightarrow (N(c),g_N)$ be a Riemannian map from a Riemannian manifold to a space form. Then $F$ is harmonic if and only if  $kerF_\ast$ is minimal.
	\end{theorem}
	\begin{proof} Since $(N,g_N)$ be a Ricci soliton then, we have
		\begin{equation*} 
			\frac{1}{2} (L_{F_\ast X} g_N)(F_\ast Y,V) + Ric (F_\ast Y,V) + \lambda g_N(F_\ast Y,V) = 0,
		\end{equation*}
		for $V \in \Gamma(rangeF_\ast)^\bot$ and $F_\ast X, F_\ast Y \in \Gamma(rangeF_\ast)$. Then from above equation, we get
		\begin{equation*}
			\frac{1}{2} \{g_N(\nabla_{F_\ast Y}^N F_\ast X, V) + g_N(\nabla_V^N F_\ast X, F_\ast Y) \} + Ric (F_\ast X, V) = 0,
		\end{equation*}
		or
		\begin{equation}\label{eqn5.7} 
			\frac{1}{2} \{g_N(\overset{N}{\nabla_{Y}^F} F_\ast X \circ F, V) + g_N(\nabla_V^N F_\ast X, F_\ast Y) \} + Ric (F_\ast X, V) = 0.
		\end{equation}
		By definition of Ricci tensor and using (\ref{eqn2.2}) in (\ref{eqn5.7}), we get
		\begin{equation*}
			\begin{array}{ll}
				\frac{1}{2} \{g_N(F_\ast (\nabla_Y^M X) + (\nabla F_\ast)(Y,X), V) + g_N(\nabla_V^N F_\ast X, F_\ast Y) \} \\+ \sum\limits_{j=r+1}^{m}g_N(R^N(F_\ast X_j, F_\ast X)V, F_\ast X_j) + \sum\limits_{k=1}^{n_1} g_N(R^N(e_k, F_\ast X)V, e_k)  = 0,
			\end{array}
		\end{equation*}
		where $\{F_\ast X_j\}_{r+1 \leq j \leq m}$ and $\{e_k\}_{1 \leq k \leq n_1}$ are orthonormal bases of $rangeF_\ast$ and $(rangeF_\ast)^\bot$, respectively. Using (\ref{eqn2.10}) in above equation, we get
		\begin{equation}\label{eqn5.8} 
			\begin{array}{ll}
				\frac{1}{2} \{g_N( (\nabla F_\ast)(Y,X), V) + g_N(\nabla_V^N F_\ast X, F_\ast Y) \} \\+ \sum\limits_{j=r+1}^{m}c \Big \{ g_N( g_N(F_\ast X, V) F_\ast X_j - g_N(F_\ast X_j, V) F_\ast X, F_\ast X_j) \Big \} \\+ \sum\limits_{k=1}^{n_1} c \Big \{ g_N( g_N(F_\ast X, V) e_k - g_N(e_k, V) F_\ast X, e_k) \Big \} = 0.
			\end{array}
		\end{equation}
		Taking trace of (\ref{eqn5.8}), we get
		\begin{equation}\label{eqn5.9} 
			\begin{array}{ll}
				\frac{1}{2} \sum\limits_{j=r+1}^{m}\underset{k=1}{\overset{n_1}{\sum}} \Big \{g_N( (\nabla F_\ast)(X_j,X_j), e_k) + g_N(\nabla_{e_k}^N F_\ast X_j, F_\ast X_j) \Big \}  = 0.
			\end{array}
		\end{equation}
		Since $\nabla^N$ is metric connection on $N$ and using (\ref{eqn2.2}) in (\ref{eqn5.9}), we get
		\begin{equation*}
			\begin{array}{ll}
				\sum\limits_{j=r+1}^{m}\underset{k=1}{\overset{n_1}{\sum}} g_N( \overset{N}{\nabla_{X_j}^F} F_\ast(X_j), e_k) = 0.
			\end{array}
		\end{equation*}
		Now using (\ref{eqn5.1}) in above equation, we obtain
		\begin{equation*}
			\begin{array}{ll}
				\underset{k=1}{\overset{n_1}{\sum}} g_N( H_2, e_k) = 0.
			\end{array}
		\end{equation*}
		Hence $H_2 = 0$. Then by Lemma \ref{lem5.1}, $F$ is harmonic if and only if $H=0$, which completes the proof.
	\end{proof}
	Eells and Sampson introduced the notion of biharmonic map in \cite{Eells}. A map $F:(M^m,g_M)\rightarrow (N^n,g_N)$ between Riemannian manifolds is biharmonic if and only if the bitension field of $F$ vanishes at each point $p\in M$, i.e.
	\begin{equation*}
		\tau_2(F)= -\Delta^F \tau(F) - trace_{g_M} R^N(dF, \tau(F)) dF= 0.
	\end{equation*}
	In other words, biharmonic map is a critical point of bienergy. Further, Jiang obtained Euler-Lagrange equations for biharmonic map in \cite{Jiang1986a}. The biharmonicity of immersions and submersions was studied in \cite{Jiang1986, Oniciuc, Akyol, Urakawa}. B. \c{S}ahin studied biharmonic Riemannian maps and obtained the following necessary and sufficient condition.
	\begin{theorem} \cite{Sahin2011a} Let $ F:(M^m,g_M) \rightarrow (N(c),g_N)$ be a Riemannian map from a Riemannian manifold to a space form. Then $F$ is biharmonic if and only if
		\begin{equation}\label{eqn5.10} 
			\begin{array}{ll}
				r trace \mathcal{S}_{(\nabla F_\ast)(.,H)} F_\ast(.) - r trace F_\ast(\nabla_{(.)} \nabla_{(.)} H)\\- (m-r) trace F_\ast(\nabla_{(.)}  {}^\ast F_\ast(\mathcal{S}_{H_2} F_\ast(.))) - (m-r) trace \mathcal{S}_{\nabla_{(.)}^{F\bot} H_2} F_\ast (.) \\- r c(m-r-1) F_\ast(H) = 0,
			\end{array}
		\end{equation}
		and
		\begin{equation}\label{eqn5.11} 
			\begin{array}{ll}
				r trace \nabla_{(.)}^{F \bot} {(\nabla F_\ast)(.,H)} + r trace	(\nabla F_\ast)(., \nabla_{(.)} H)\\ + (m-r) trace (\nabla F_\ast) (., {}^\ast F_\ast(\mathcal{S}_{H_2} F_\ast(.))) - (m-r) \Delta^{R^\bot} H_2 \\- (m-r)^2 c H_2 = 0,
			\end{array}	
		\end{equation}
		where $dim(kerF_\ast) = r$ and $dim(kerF_\ast)^\bot = m-r$.
	\end{theorem}
	\begin{theorem} Let $(N,g_N,F_\ast X,\lambda)$ be a Ricci soliton with potential vector field $F_\ast X \in \Gamma(rangeF_\ast)$ and $F:(M,g_M) \rightarrow (N(c),g_N)$ be a Riemannian map from a Riemannian manifold to a space form. Then $F$ is biharmonic if and only if $kerF_\ast$ is minimal.
	\end{theorem}
	\begin{proof} We see in Theorem \ref{thm5.2}, $H_2 =0$ then from (\ref{eqn5.10}) and (\ref{eqn5.11}), $F$ is biharmonic if and only if $H=0$, which completes the proof.
	\end{proof}
	\section{Example}\label{sec6}
	\begin{example}
		Let $M= \{ (x_1, x_2, x_3) \in \Bbb{R}^{3} : x_1 \neq 0, x_2 \neq 0, x_3 \neq 0\}$ be a 3-dimensional Riemannian manifold with Riemannian metric $g_M$ on $M$ given by $g_M= e^{2x_3}dx_1^2 + e^{2x_3}dx_2^2 + e^{2x_3} dx_3^2$. Let $N=\{ (y_1, y_2) \in \Bbb{R}^{2} \}$ be a Riemannian manifold with Riemannian metric $g_N$ on $N$ given by $g_N= e^{2x_3}dy_1^2 + dy_2^2$. Consider a map $F : (M,g_M) \rightarrow (N,g_N) $ defined by
		\begin{equation*}
			F(x_1,x_2,x_3)= \Big( \frac{x_1 + x_2 + x_3}{\sqrt {3}}, 0 \Big).
		\end{equation*}
		By direct computations
		\begin{equation*}	
			kerF_\ast = Span \Big\{ U_1 = -e_1+e_2, U_2 = -e_1+e_3 \Big\}
		\end{equation*}
		and 
		\begin{equation*}
			(kerF_\ast)^\bot = Span \Big\{ X=\frac{e_1+e_2+e_3}{\sqrt{3}} \Big\},
		\end{equation*}
		where $\Big\{ e_1 = e^{-x_3}\frac{\partial}{\partial x_1}, e_2 = e^{-x_3}\frac{\partial}{\partial x_2}, e_3 = e^{-x_3}\frac{\partial}{\partial x_3} \Big\}$, $\Big\{ e_1' = e^{-x_3}\frac{\partial}{\partial y_1}, e_2' = \frac{\partial}{\partial y_2} \Big\}$ are bases on  $T_pM$ and $T_{F (p)}N$ respectively, for all $ p\in M$. By direct computations, we can see that $F_\ast (X) = e_1'$ and $ g_M(X,X)= g_N(F_\ast X, F_\ast X)$ for $X \in \Gamma(kerF_\ast)^\bot.$ Thus $F$ is a Riemannian map with $rangeF_\ast= Span \Big\{ F_\ast X = e_1'\Big\}$ and $(rangeF_\ast)^\bot= Span \Big\{ e_2'\Big\}.$ Now, we will show that base manifold $N$ admits a Ricci soliton, i.e.
		\begin{equation}\label{eqn6.1} 
			\frac{1}{2} (L_{Z_1} g_N)(X_1,Y_1) + Ric (X_1,Y_1) + \lambda g_N(X_1,Y_1) = 0,
		\end{equation}
		for any $X_1, Y_1, Z_1 \in \Gamma(TN)$. Now,
		\begin{equation}\label{eqn6.2} 
			\frac{1}{2} (L_{Z_1} g_N)(X_1,Y_1) = \frac{1}{2} \Big \{ g_N(\nabla_{X_1}^N Z_1,Y_1) + g_N(\nabla_{Y_1}^N Z_1,X_1) \Big \}.
		\end{equation}
		Since dimension of $rangeF_\ast$ and $(rangeF_\ast)^\bot$ is one therefore we can decompose $X_1, Y_1$ and $Z_1$ such that $X_1= a_1e_1'+a_2e_2'$, $Y_1=a_3e_1'+a_4e_2'$ and $Z_1= a_5e_1'+a_6e_2'$, where $e_1'$ and $e_2'$ denote for components on $rangeF_\ast$ and $(rangeF_\ast)^\bot$, respectively and $\{a_i\}_{1 \leq i \leq 6} \in \Bbb{R}$ are some scalars. Then from (\ref{eqn6.2}), we get
		\begin{equation*}
			\begin{array}{ll}
				\frac{1}{2} (L_{Z_1} g_N)(X_1,Y_1)=& \frac{1}{2} \Big \{ g_N(\nabla_{a_1e_1'+a_2e_2'}^N a_5e_1'+a_6e_2',a_3e_1'+a_4e_2') \\&+ g_N(\nabla_{a_3e_1'+a_4e_2'}^N a_5e_1'+a_6e_2',a_1e_1'+a_2e_2') \Big \}.
			\end{array}
		\end{equation*}
		Since $\nabla^N$ is metric connection then from above equation, we get
		\begin{equation}\label{eqn6.3} 
			\begin{array}{ll}
				\frac{1}{2} (L_{Z_1} g_N)(X_1,Y_1)= &\frac{1}{2}\Big \{  2 a_1a_3a_6g_N(  \nabla_{e_1'}^N e_2', e_1') + 2 a_2a_4a_5g_N(  \nabla_{e_2'}^N e_1', e_2') \\&+a_2a_3a_6 g_N( \nabla_{e_2'}^N e_2', e_1') + a_1 a_4 a_5 g_N(\nabla_{e_1'}^N e_1', e_2') \\&+ a_1 a_4 a_6 g_N(\nabla_{e_2'}^N e_2', e_1') +a_2a_3a_5 g_N( \nabla_{e_1'}^N e_1', e_2') \Big \}.
			\end{array}
		\end{equation}
		Since $\nabla_{e_1'}^N e_1' =0$, $\nabla_{e_1'}^N e_2' =0$, $\nabla_{e_2'}^N e_1' =0$ and $\nabla_{e_2'}^N e_2' =0$, by (\ref{eqn6.3}), we get
		\begin{equation}\label{eqn6.5} 
			\frac{1}{2} (L_{Z_1} g_N)(X_1,Y_1) = 0.
		\end{equation}
		Also,
		\begin{equation}\label{eqn6.6} 
			g_N(X_1,Y_1)= g_N(a_1e_1'+a_2e_2',a_3e_1'+a_4e_2') = (a_1a_3 + a_2a_4),
		\end{equation}
		and
		\begin{equation*}
			Ric(X_1,Y_1) = Ric(a_1e_1'+a_2e_2',a_3e_1'+a_4e_2'),
		\end{equation*}
		which implies
		\begin{equation}\label{eqn6.7} 
			Ric(X_1,Y_1)= a_1a_3Ric(e_1',e_1') + (a_1a_4+a_2a_3) Ric(e_1',e_2') +a_2a_4 Ric(e_2',e_2').
		\end{equation}
		By (\ref{eqn3.10}) and (\ref{eqn3.11}), we get
		\begin{equation}\label{eqn6.8} 
			\begin{array}{ll}
				Ric(e_1', e_1')=&  Ric^{rangeF_\ast}(e_1', e_1') -  g_N(\mathcal{S}_{\nabla_{e_2'}^{F \bot} e_2'} e_1', e_1') \\&+  g_N(\nabla_{e_2'}^N \mathcal{S}_{e_2'} e_1' , e_1') - g_N(\mathcal{S}_{e_2'}  e_1', \mathcal{S}_{e_2'} e_1')  \\&- g_N(\nabla_{e_2'}^N  e_1', \mathcal{S}_{e_2'} e_1'),
			\end{array}
		\end{equation}
		and
		\begin{equation}\label{eqn6.9} 
			\begin{array}{ll}
				Ric(e_2', e_2')=&  Ric^{(rangeF_\ast)^\bot}(e_2', e_2') - g_N(\mathcal{S}_{\nabla_{e_2'}^{F \bot} e_2'} e_1', e_1') \\&+  \nabla_{e_2'}^N ( g_N( \mathcal{S}_{e_2'} e_1' , e_1')) - g_N(\mathcal{S}_{e_2'}  e_1', \mathcal{S}_{e_2'} e_1')  \\&- 2g_N(\nabla_{e_2'}^N  e_1', \mathcal{S}_{e_2'} e_1').
			\end{array}
		\end{equation}
		By (\ref{eqn2.8}) and (\ref{eqn3.9}), we get
		\begin{align}\label{eqn6.10} 
				Ric(e_1', e_2')&= g_N(R^N(e_1', e_1') e_2', e_1') + g_N(R^N(e_2', e_1') e_2', e_2') \nonumber \\& = g_N(\nabla_{e_2'}^N \nabla_{e_1'}^N e_2' - \nabla_{e_1'}^N \nabla_{e_2'}^N e_1' - \nabla_{[e_2', e_1']}^N e_2', e_2')=0.
		\end{align}
		Using (\ref{eqn6.8}), (\ref{eqn6.9}) and (\ref{eqn6.10}) in (\ref{eqn6.7}), we get
		\begin{equation}\label{eqn6.11} 
			\begin{array}{ll}
				Ric(X_1,Y_1)=&  (a_1a_3)Ric^{rangeF_\ast}(e_1', e_1') +(a_1a_3)  g_N(\nabla_{e_2'}^N \mathcal{S}_{e_2'} e_1' , e_1') \\&- (a_1a_3+a_2a_4) g_N(\mathcal{S}_{e_2'}  e_1', \mathcal{S}_{e_2'} e_1')   \\&+ (a_2a_4)Ric^{(rangeF_\ast)^\bot}(e_2', e_2') \\&+ (a_2a_4) \nabla_{e_2'}^N ( g_N( \mathcal{S}_{e_2'} e_1' , e_1')).
			\end{array}
		\end{equation}
		Since dimension of $rangeF_\ast$ and $(rangeF_\ast)^\bot$ is one therefore  $Ric^{rangeF_\ast}(e_1', e_1')= 0$ and $Ric^{(rangeF_\ast)^\bot}(e_2', e_2') = 0$. Also since  $\mathcal{S}_{e_2'} e_1' \in \Gamma(rangeF_\ast)$, we can write $\mathcal{S}_{e_2'} e_1' = a_7 e_1'$ for some scalar $a_7 \in \Bbb{R}$. Then by substituting these values in (\ref{eqn6.11}), we get
		\begin{equation*}
			\begin{array}{ll}
				Ric(X_1,Y_1)=& (a_1a_3)  g_N(\nabla_{e_2'}^N a_7 e_1' , e_1') \\&- (a_1a_3+a_2a_4) g_N(a_7  e_1', a_7 e_1') \\&+ (a_2a_4) \nabla_{e_2'}^N ( g_N( a_7 e_1' , e_1')),
			\end{array}
		\end{equation*}
		which implies
		\begin{equation}\label{eqn6.12} 
			\begin{array}{ll}
				Ric(X_1,Y_1)= - (a_1a_3+a_2a_4) a_7^2.
			\end{array}
		\end{equation}
		Now, using (\ref{eqn6.5}), (\ref{eqn6.6}) and (\ref{eqn6.12}) in (\ref{eqn6.1}), we obtain that $N$ admits Ricci soliton for
		\begin{equation*}
				\lambda = a_7^2.
		\end{equation*}	
	\end{example}
	
	\noindent \textbf{Acknowledgment:} 
	We appreciate and thank to the referees for their questions, comments and suggestions to improve the quality of the paper. We are very grateful to them.

\end{document}